\documentclass[preprint,12pt]{elsarticle}

\usepackage{amssymb}
\usepackage{amsthm}
\usepackage{amsmath}

\usepackage{color}

\journal{Journal of Computational and Applied Mathematics}

\begin{document}

\begin{frontmatter}



\title{Simultaneous {Single-Step} One-Shot Optimization with Unsteady PDEs}


\author[TUKL]{Stefanie G\"unther\corref{Guenther}}
\cortext[Guenther]{Corresponding author}
\ead{stefanie.guenther@scicomp.uni-kl.de}

\author[TUKL]{Nicolas R. Gauger}
\address[TUKL]{TU Kaiserslautern, Chair for Scientific Computing, Paul-Ehrlich-Stra{\ss}e 34, 67663 Kaiserslautern, Germany}

\author[MIT]{Qiqi Wang}
\address[MIT]{Massachusetts Institute of Technology, Department of Aeronautics and Astronautics, 77 Massachusetts Avenue, Cambridge, MA 02139, USA}

\begin{abstract}
The {single-step} one-shot method has proven to be very efficient for PDE-constrained optimization where the partial differential equation (PDE) is solved by an iterative fixed point solver. In this approach, the simulation and optimization tasks are performed simultaneously in a single iteration. If the PDE is unsteady, finding an appropriate fixed point iteration is non-trivial. In this paper, we provide a framework that makes the {single-step} one-shot method applicable for unsteady PDEs that are solved by classical time-marching schemes. The one-shot method is applied to an optimal control problem with unsteady incompressible Navier-Stokes equations that are solved by an industry standard simulation code. With the Van-der-Pol oscillator as a generic model problem, the modified simulation scheme is further improved using adaptive time scales. Finally, numerical results for the advection-diffusion equation are presented.

\end{abstract}


\end{frontmatter}



\section{Simultaneous PDE-constrained optimization}
\label{sec:intro}

Many engineering and science problems can be described by partial differential equations (PDEs). Instead of solving them analytically, the rapid progress in computer technologies made it possible
to compute approximate solutions numerically. In Computational Fluid Dynamics (CFD), advanced and mature simulation tools have been developed, refined and validated over years. Nowadays, high fidelity simulation tools are available that compute accurate approximations for a variety of complex fluid flow configurations \cite{nielson2014fun3d, kroll2014tau}.

The simulation task refers to the process of solving the PDE numerically for some state variables (as for example the velocity, density or temperature of a system) while some appropriate data is given (such as geometry, material coefficients, boundary or initial conditions). In contrast, the task of optimization is to adjust some of this data in such a way, that the state variables exhibit a desired behavior determined by an objective function $J$. Let $y$ be the vector of state variables of the system and let $u$ describe the data that can be adjusted, the so-called design variables. The optimization problem under consideration then reads
\begin{align}\label{general_min_problem}
	\min_{y,u} \, J(y,u) \quad \text{s.t.} \quad c(y,u) = 0 
\end{align}
where $c$ represents a system of PDEs including boundary and/or initial conditions. Depending on the choice of the state and design variables, the optimization problem can represent an optimal shape design problem, an inverse design problem for parameter estimation or an optimal control problem. It has a wide range of application scenarios as for example in aerodynamic shape design, where one aims to find an airfoil shape that minimizes its drag coefficient while the fluid flow satisfies the Navier-Stokes equations, or applications in other disciplines like geophysics, medical imaging and atmospheric science. 

If evaluating the objective function is rather time consuming, the optimization problem is typically solved using {gradient-based methods \cite{nocedal2006numerical,lions1971optimal}}. These methods iteratively perform design updates utilizing the sensitivity of the objective function to the design variables. If the dimension of the design space is rather large, the adjoint method is preferred since the cost for computing the sensitivities is then independent of the dimension of the design space \cite{pironneau1974optimum, nadarajah2007optimum}. In this approach, an adjoint PDE system is derived from the Lagrangian function associated with \eqref{general_min_problem}. The solution of the adjoint system can be used to determine the desired sensitivities efficiently. Two main approaches are distinguished namely the discrete adjoint approach, where a discrete adjoint system is derived from the discretized PDEs, and the continuous adjoint approach that derives the adjoint system from the continuous PDEs and discretizes afterwards. The discrete adjoint solution can be generated automatically using Automatic Differentiation (AD) \cite{griewank2008evaluating}.

In most application scenarios it is reasonable to assume, that every design variable $u$ uniquely determines a state $y(u)$ that satisfies the state equation $c(y(u),u)=0$. Under this assumption it is a common approach to eliminate the state variable from the optimization and focus on the unconstrained minimization problem 
\begin{align}\label{reduced_optim_problem}
	\min_u \, J(y(u),u).
\end{align}
 In this so called \textit{reduced space approach} a gradient-based optimization strategy can be applied to the reduced objective function \eqref{reduced_optim_problem} that depends on the design $u$ solely, while the PDE-constraint is treated implicitly by recovering $c(y(u),u)=0$ after each design change. Methods of this type are also referred to as \textit{black-box approach} or \textit{Nested Analysis and Design (NAND)} \cite{akccelik2006parallel}. The main drawback of reduced space methods is a direct consequence of the implicit treatment of the constraint: After each design change, the state variable has to be recomputed such that it satisfies the PDE. This means, that a full numerical simulation has to be performed after each design change. Despite of the rapid growth in computer capacities, simulating nonlinear PDE-systems still can take hours or even weeks on state-of-the-art supercomputers. This makes the reduced space approach unaffordable in many sophisticated application scenarios. 

On the other extreme, so called \textit{full space methods} are a popular alternative \cite{biegler2003large}. Instead of reducing the optimization space by recovering the PDE-solution after each design change, the optimization problem is solved in the full space for the (discretized) state and design variables. In this approach, the optimality conditions for the constrained minimization problem \eqref{general_min_problem} are solved simultaneously for the state, the adjoint and the design variable in a SQP-like fashion. Because the simulation is directly integrated in the optimization process, these methods are often  called \textit{Simultaneous Analysis and Design} (SAND), \textit{all-at-once approach}, or \textit{one-shot approach} {\cite{ziems2011adaptive,biros2005paralleli,biros2005parallelii, abbeloos2011nested, Taasan1991oneshot,gauger2009singlestep,Hazra2004pseudotime}}. It has been observed numerically - at least for steady state PDEs - that the full space methods can outperform the reduced space methods by about one order of magnitude measured in iteration counts and runtime \cite{akccelik2006parallel, gauger2012automated}. However, the major drawback of many full space optimization methods is, that they require the computation of additional Jacobians and Hessians which are not necessarily part of the PDE simulation tool. In fact, many PDE solvers only approximate Jacobians of the discretized PDE residuals due to implementation or computational effort which makes it necessary to rewrite and enhance the PDE solver for optimization. 

{In this paper, we consider a full space optimization method that overcomes this difficulty. We follow the single-step one-shot approach as proposed in \cite{gauger2012automated,gauger2009singlestep} that is specially tailored for PDE-constrained optimization where it is impossible to form or factor the Jacobian of the constraint. Instead, it is assumed that the user is provided with an iterative fixed point algorithm that computes a discrete numerical approximation to the PDE solution in a black box fashion.} In the {considered} one-shot approach, these iterations are enriched by an iteration for the adjoint as well as the design variable, {so that in each optimization iteration, only one step of the PDE solver and the adjoint solver is executed. The iterative adjoint solver} can be computed efficiently with the reverse mode of AD applied to the PDE solver and evaluating the objective function. The design step is based on the gradient of the reduced space objective function. Provided that a certain preconditioner for the design update is used, the {single-step} one-shot method is proven to converge to an optimal point of the minimization problem \cite{griewank2010properties}. Since the iterations of the PDE fixed point solver are used in a black-box manner, the {optimization method} leverages and retains the software investment that has been made in developing the PDE solver. Section \ref{sec:oneshot} shortly recalls the main aspects of the {considered one-shot approach}. 

Application of the {single-step} one-shot method to optimization with steady state PDEs is {straightforward} in terms of the fixed point solver: It is a common and well established approach for solving steady state PDEs to apply the so-called pseudo-time-stepping method. In this approach, the PDE is interpreted as a steady state of a dynamical system and solved numerically by an explicit (pseudo-)time-stepping method \cite{schulz2004aerodynamic, Hazra2004pseudotime}. Its strong relation to general iterative methods made it possible to apply the proposed one-shot method to various optimization tasks with steady state PDEs especially in computational fluid dynamics \cite{gauger2012automated, bosse2014optimal,gauger2009singlestep,griewank2011reduced}. 

However, if the PDE is fully unsteady, finding an appropriate fixed point iteration that fits in the {single-step} one-shot framework raises complexity. Existing simulation tools for unsteady PDEs typically apply an implicit time marching scheme. The resulting implicit equations are solved one after another forward in time utilizing an iterative solver at each time step as proposed in the famous dual time-stepping approach by Jameson \cite{jameson1991dualtime}. In order to prepare for {single-step} one-shot optimization, a fixed point iterator for the unsteady PDE is derived from such a method in Section \ref{sec:unsteadyOneshot} by reducing the number of {iterations at each time step. Since the iteration steps themselves are not changed, the effort that has been spent} for developing the PDE-solver is preserved within the new scheme. 

In Section \ref{sec:unsteadyRANS}, the {proposed} one-shot method is applied to an optimal control problem with unsteady incompressible Navier-Stokes equations that are solved by an industry standard simulation code. The test case under consideration is an active flow control problem of a cylinder in unsteady flow. Actuation slits are installed on the cylinder surface in order to reduce vorticity in the wake.

The modified time marching scheme is further improved using adaptive time scales in Section \ref{sec:improving}. Numerical tests on model problems have shown, that the number of iterations needed to solve the unsteady PDE can thereby be reduced drastically \cite{guenther2014extension}. The adaptive time scale approach is further investigated in Section \ref{sec:numresult} for the advection-diffusion equation with periodic boundary condition.

\section{The {single-step} one-shot method}
\label{sec:oneshot}

Let $c(y,u) = 0$ with $c\colon Y\times U \to Y$ represent a system of PDEs with state vector $y$ and a set of design variables $u\in U$. In a discretized PDE setting, it is reasonable to assume that $Y$ and $U$ are finite dimensional Hilbert spaces with $\dim Y=m$ and $\dim U=n$ which allows us to associate their elements with the corresponding coordinate vectors in $\mathbb{R}^m$ and $\mathbb{R}^n$, respectively, and write duals as transposed vectors. For an objective function $J\colon Y\times U\to \mathbb{R}$ we want to solve the PDE-constrained optimization problem
\begin{align}\label{minproblem}
	\min_{y,u} \, J(y,u) \quad \text{s. t.} \quad c(y,u) = 0. 
\end{align}

The {proposed} one-shot method is tailored for optimization problems where the PDE-constraint is solved by an iterative fixed point solver that computes a feasible state variable given appropriate design data. We therefore assume, that an iteration function $H\colon Y\times U \to Y$ is available which - for any given design $u^*\in U$ - iteratively updates the state variable with
\begin{align}\label{steady_primaliteration}
y_{k+1}=H(y_k,u^*)
\end{align}
such that the limit $y^* = \lim_{k\to \infty} y_k \in Y$ exists and satisfies $c(y^*,u^*)=0$. Convergence of the above iteration is assured if
\begin{align}
\left\|\frac{\partial H}{\partial y} \right\| \leq \rho  <1 
\end{align}
holds for all points of interest, i.e. the iteration function $H$ is contractive with respect to the state variable according to Banach's fixed point theorem \cite{barner1995analysis}. If $\rho$ is close to $1$, the simulation code is a rather slowly converging fixed point iteration as it is typically the case for applications in CFD. 

Assuming that the limit $y^*$ is a fixed point of the iteration function $H$ if and only if $c(y^*,u^*)=0$, we can reformulate the constraint of the optimization problem in terms of the fixed point equation $y^*=H(y^*,u^*)$ and focus on the following minimization problem
\begin{align}\label{minproblem_fixedpoint}
	\min_{y,u} \, J(y,u) \quad \text{s. t.} \quad y=H(y,u).
\end{align}
We define the associated Lagrangian function
\begin{align}
	L(y,\bar y,u):=J(y,u)+(H(y,u)-y)^T\bar y
\end{align}
with Lagrange multiplier $\bar y\in Y^*$ which corresponds to the so-called adjoint variable. Any local optimizer of \eqref{minproblem_fixedpoint} is a saddle-point of $L$  
 leading to the following necessary optimality conditions, the so-called Karush-Kuhn-Tucker (KKT) conditions \cite{nocedal2006numerical}:
\begin{align} 
 y&=H(y,u)  & \textit{state equation} \label{NOC1} \\
\bar y &= J_y(y,u)^T +  H_y(y, u)^T\bar y & \textit{adjoint equation}  \label{NOC2}\\ 
 0 &= J_u(y,u)^T + H_u(y,u)^T\bar y & \textit{design equation} \label{NOC3}
\end{align}
where subscripts denote partial derivatives. The state equation, also referred to as the primal equation, can be solved with the fixed point iteration \eqref{steady_primaliteration} by construction. The adjoint equation is linear in the adjoint variable $\bar y$ and involves the transpose of the Jacobian of the primal fixed point iterator with respect to the state. An iteration for solving the adjoint equation can be obtained by applying AD to $H$ and evaluating $J$ \cite{griewank2008evaluating}. This approach automatically generates an iteration of the following type
\begin{align}\label{steady_adjointiteration}
	\bar y_{k+1} =  J_y(y,u)^T + H_y(y, u)^T\bar y_k.
\end{align}
Since $H$ is contractive, this iteration converges to a solution of \eqref{NOC2} for any given design and corresponding state that satisfy the primal equation. The design equation refers to stationarity of the objective function with respect to design changes. For feasible state and adjoint variables, it corresponds to the total derivative of reduced space objective function in \eqref{reduced_optim_problem} and is often called the reduced gradient. The reduced gradient is therefore utilized for updating the design variable during the optimization procedure. Again, all partial derivatives can be evaluated efficiently with the use of AD.

Instead of solving the primal equation with \eqref{steady_primaliteration} first, then iterating for an adjoint solution with \eqref{steady_adjointiteration}, followed by an update of the design in the direction of the reduced gradient, the main idea of the {single-step} one-shot method is to solve the set of KKT-equations simultaneously in one single iteration \cite{griewank2010properties,griewank2011reduced}:
\begin{eqnarray} \label{oneshot_iteration}
  \begin{bmatrix} y_{k+1} \\ \bar y_{k+1} \\ u_{k+1} \end{bmatrix} = \begin{bmatrix} H(y_k,u_k) \\ J_y(y_k,u_k)^T + H_y(y_k, u_k)^T\bar y_k \\ u_k - B_k^{-1} \left(J_u(y_k,u_k)^T + H_u(y_k,  u_k)^T\bar y_k,\right) \end{bmatrix}\, .
\end{eqnarray}
In this approach, the adjoint iteration and a preconditioned reduced gradient step for the design variable are integrated into the primal iteration. In order to ensure convergence of the {single-step} one-shot iteration, the preconditioner $B_k$ has to be chosen such that the coupled iteration is contractive. It is suggested in \cite{griewank2010properties,griewank2011reduced} to look for descent on the doubly augmented Lagrangian function
\begin{multline}\label{def:augmented_Lagrangian}
  L^a(y,\bar y,u) :=  \frac{\alpha}{2}\left\|H(y,u)-y\right\|^2 + \frac{\beta}{2}\left\|J^T_y(y,u) + (H_y(y,u)-I)^T\bar y\right\|^2 \\ + L(y,\bar y, u)
\end{multline}
where weighted residuals of the state and the adjoint equations are added to the Lagrangian function with weights $\alpha, \beta >0$. It is proven ibidem, that a suitable preconditioner approximates the Hessian $L_{uu}^a$ of the augmented Lagrangian with respect to the design $u$. {Numerically, $B_k$ can therefore be approximated using secant updates on the gradient of the augmented Lagrangian $\nabla_uL^a$, as for example applying BFGS-updates in each iteration \cite{nocedal2006numerical}. Computation of the gradient can be automated with the use of AD, where \textit{forward over reverse} differentiation provides routines for computing second derivatives \cite{naumann2012art}.}

The above {single-step} one-shot method updates the primal, the adjoint and the design variable simultaneously in a Jacobi-like fashion. In \cite{bosse2014oneshot} variants of the method are surveyed which update the variables in a Seidel-type iteration where the new variables are used immediately. The design space preconditioner should then approximate the reduced space Hessian of the objective function. Furthermore, multi-step one-shot methods are analysed where not only one but several steps of the primal and the adjoint iteration are performed before the design is updated.

The single-step one-shot method has been applied to various optimization problems where the underlying PDE is steady \cite{gauger2012automated,griewank2011reduced,bosse2014optimal}. In that case, the fixed point iterator $H$ arises naturally in common simulation tools. It can be for example one step of an explicit pseudo-time-stepping scheme which is typically state of the art for solving steady state PDEs in industrial CFD applications \cite{nielson2014fun3d, kroll2014tau}. {Numerical applications have proven}, that the cost for an one-shot optimization is only a small multiple of the cost of a single simulation of the underlying PDE - a property which is called {\it bounded retardation}. The factor typically varies between $2$ and $8$. {Direct comparison with a classical reduced space BFGS optimization showed, that the overall runtime for the one-shot optimization is about one order of magnitude lower than for the reduced space approach \cite{oetzkaya2014oneshot,bosse2014optimal}.}

\section{{Single-step} one-shot optimization with unsteady PDEs}\label{sec:unsteadyOneshot}

With the rapid increase in computational capacities, numerical simulation codes are no longer restricted to steady state solutions but perform accurate high fidelity simulations of unsteady turbulent flow. Common simulation methods discretize the transient term of the PDE in time applying a time-stepping method (e.g. Runge-Kutta) while implicit methods are often preferred due to good stability properties {\cite{lambert1991numerical}}. The resulting implicit equations are then solved iteratively as for example in dual time-stepping methods with Picard- or Newton-iterations in each time step, often in combination with multigrid methods and Krylov subspace techniques {\cite{jameson1991dualtime,Breuer1993adualtime}}. In this section we present a framework for {single-step} one-shot optimization with unsteady PDEs utilizing these time-marching schemes.

\subsection{Problem statement}

For time-dependent PDEs the state variable varies with time and thus is a function $y\colon \mathbb{R} \to Y$. The objective function to be minimized is typically given by a time averaged quantity
\begin{align}
	J(y,u) := \frac 1 T \int_0^T \, \hat J(y(t),u) \, \mathrm{d}t
\end{align}
where $\hat J \colon Y \times U \to \mathbb{R}$ represents some time-dependent quantity as for example drag or lift coefficient of an airfoil in unsteady flow. The optimization problem with unsteady PDE-constraints then reads
\begin{align} \label{OPunsteady}
\min_{y,u} \, \frac 1 T \int_0^T& \, \hat J(y(t),u) \, \mathrm{d}t \quad \text{ subject to} \\
	\frac{\partial y(t)}{\partial t} &= f(y(t),u)   \qquad \forall\, t\in [0,T]  \\
	y(0)&=y^0 
\end{align}
where the right hand side  $f\colon Y\times U\to Y$ corresponds to {spatial} derivative operators and boundary terms of the unsteady PDE and $y^0 \in Y$ is some appropriate initial data.

\subsection{Fixed point iteration for unsteady PDE-constraints}\label{subsec:fixedpointunsteady}

Numerical methods for solving unsteady PDEs discretize the time domain into a finite number of time steps with $t_0=0 <t_1 <\dots < t_N=T$ and advance the solution forward in time. The transient term is typically discretized by an implicit scheme due to stability reasons which results in a (nonlinear) implicit residuum equation at each time step:
\begin{align}\label{unsteady_residuum_equations}
	R(y^i,y^{i-1},y^{i-2},\dots,u) \text{{$=$}}0 \quad \forall \, i = 1, \dots, N 
\end{align}
where $y^{i}\approx y(t_{i})$ approximates the solution at time $t_{i}$. 
Depending on the order of the implicit time-stepping approximation, the residuum equation for a certain time step approximation $y^i$ contains approximations to the solution at one or more previous time steps $y^{i-1}, y^{i-2}$ etc. For notational reasons, we choose the first order Backward Euler discretization and drop previous states except $y^{i-1}$. Application of the proposed method to higher order schemes is {straightforward}. In case of the Backward Euler method, the residuum equations are 
\begin{align}\label{unsteady_residuum}
	R(y^{i},y^{i-1},u) := \frac{y^{i} - y^{i-1}}{t_{i} - t_{i-1}}  - f_h(y^{i},u) \text{{$=$}}0 \quad \forall \, i=1,\dots,N
\end{align}
where $f_h$ represents a spatial discretization of the right hand side of the PDE. 

The set of nonlinear residuum equations can be solved one after another marching forward in time. Typically, iterative methods are used to converge to a pseudo-steady state solution at each time step:
\begin{align}\label{unsteady_primal_classic}
	\text{for } i=1,\dots, N: &\notag \\
 	\text{ iterate } \,  &y^{i}_{k+1}=G^i(y^{i}_k, y^{i-1}, u) \,\, \overset{k\to\infty}{\longrightarrow} \,y^{i}
\end{align}
with $y_0^i:= y^{i-1}$ where the iterator $G^i$ is designed such that the converged pseudo-steady states $y^{i}$ satisfy the residuum equations \eqref{unsteady_residuum}. We therefore assume, that $G^i$ is contractive with respect to $y^i$, i.e.
\begin{align}\label{unsteady_contractivity_of_G}
	 \left\|\frac{\partial G^i(y^{i},y^{i-1},u)}{\partial y^{i}} \right\| \leq \rho^i <1 \quad \forall \, i=1,\dots,N
\end{align}
for all points of interest, which ensures convergence of the above iterations. The converged pseudo-steady states are fixed points of $G^i$ such that $y^{i}=G^i(y^{i}, y^{i-1}, u)$ holds if and only if $R(y^i,y^{i-1},u) = 0$. 

In order to extend from simulation to {single-step} one-shot optimization, where one incorporates design updates already during the primal flow computation, the time-marching scheme \eqref{unsteady_primal_classic} is modified in such a way, that the residuum equations at each time step are solved inexactly. Instead, an outer loop is performed that updates the state at all time steps:
\begin{align}\label{unsteady_primal_oneshot}
\text{iterate } k=0,1\dots \, : \notag \\
y^{i}_{k+1} &= G^i(y^{i}_k, y^{i-1}_{k+1},u) \quad \text{for} \quad i=1,\dots,N
\end{align}
with $y^0_k:=y^0\, \,\forall \, k\in \mathbb{N}$. In contrast to \eqref{unsteady_primal_classic}, where fixed point iterations are performed at each time step to reach the converged states $y^i$ one after another, in the one-shot framework \eqref{unsteady_primal_oneshot} a complete trajectory of the unsteady solution is updated within one iteration. Interpreting the time-dependent state variable as a discrete vector from the product space $\boldsymbol y=(y^1,\dots,y^N)\in Y^N:=Y\times\dots\times Y$  we can write \eqref{unsteady_primal_oneshot} in terms of an update function 
\begin{align}
\boldsymbol y_{k+1}=H(\boldsymbol y_k,u)
\end{align}
where $H\colon Y^N\times U \to Y^N$ performs the update formulas \eqref{unsteady_primal_oneshot} and is defined as
 \begin{eqnarray}\label{unsteady_define_H}
 H(\boldsymbol y,u)
 := \begin{pmatrix} G^1(y^1,y^0,u) \\ G^2(y^2,G^1(y^1,y^0,u),u) \\ \vdots \\ G^N(y^N,G^{N-1}(y^{N-1},G^{N-\text{{2}}}(y^{N-2}, \dots , G^1(y^1,y^0,u),u)\dots ,u),u) \end{pmatrix}
\end{eqnarray}

It is shown in \cite{guenther2014extension}, that the recursive iteration function $H$ is contractive with respect to $\boldsymbol y$, i.e.
\begin{align}\label{unsteady_contractivity_of_H}
 \left\|\frac{\partial H(\boldsymbol y,u)}{\partial \boldsymbol y} \right\| \leq \rho <1  
\end{align}
for all points of interest with $\rho := \max_i \rho^i$. This ensures convergence of the modified time-marching scheme \eqref{unsteady_primal_oneshot} to the fixed point $\boldsymbol y=H(\boldsymbol y,u)$ where $y^{i}=G^i(y^{i}, y^{i-1}, u)$ holds for all $i=1,\dots, N$. By construction of $H$, the fixed point satisfies the residuum equations \eqref{unsteady_residuum} and is therefore a numerical solution of the unsteady PDE.

\subsection{{Single-step} one-shot optimization steps}

The {proposed} one-shot method aims at solving the following discrete optimization problem 
\begin{align}\label{unsteady_discrete_OP}
	 \min_{\boldsymbol y,u} \, J_N(\boldsymbol y,u) \quad  \text{s.t.} \quad \boldsymbol y = H(\boldsymbol y, u) 
\end{align}
where $J_N$ is an approximation of the time averaging objective function as for example
\begin{align}
	J_N(\boldsymbol y,u) := \frac 1T \sum_{i=1}^N \Delta t_i \hat J (y^i,u) \approx \frac 1T \int_0^T \hat J(y(t),u)\, \mathrm{d}t 
\end{align}
with $\Delta t_i := t_i - t_{i-1}$. We define the corresponding Lagrangian function 
\begin{align*}
	L(\text{{$\boldsymbol y,\boldsymbol{\bar y}$}}, u) := J_N(\boldsymbol y,u) + (H(\boldsymbol y,u) - \boldsymbol y)^T\boldsymbol {\bar y}
\end{align*}
where the adjoint variable is an element from the product space $\boldsymbol{\bar y }\in (Y^*)^N$. The necessary optimality conditions for the discrete optimization problem \eqref{unsteady_discrete_OP} yield the state, the adjoint and the design equation, each of which is analyzed in the sequel:
\begin{enumerate}
\item \textit{State equation}: Differentiation of $L$ with respect to the adjoint variables yields the unsteady residuum equations at each time step. As shown in the previous subsection, these equations can be solved simultaneously in one iteration that updates the state on the entire time domain:
\begin{align}
 \text{iterate } k=0,1\dots \, : \notag  \\
 \boldsymbol y_{k+1}&=H(\boldsymbol y_k,u). 
 \end{align}

\item \textit{Adjoint equation}: From $\nabla_{y^i}L=0 $ for all $ i=1,\dots,N$, the adjoint equation is derived:
\begin{align*}
	\boldsymbol{\bar y} = \nabla_{\boldsymbol y} J_N(\boldsymbol y, u) + \left(\frac{\partial H}{\partial \boldsymbol y}\right)^T\boldsymbol{\bar y}.
\end{align*}
Since the transpose of the constraint Jacobian is an upper triangular matrix, the adjoint flow is backwards in time while $\bar y^N:=0$ is imposed at the last time step instead of initial conditions. In classical reduced space methods, the adjoint equations are solved iteratively one after another backwards in time. However, since we want to integrate the iteration into an one-shot optimization process, we rather solve the equations simultaneously for all time steps:
\begin{align}\label{unsteady_simult_adjoint}
 \text{iterate } k=0,1\dots \, : \notag \\
 \boldsymbol{\bar y}_{k+1} &=  \nabla_{\boldsymbol y} J_N(\boldsymbol y,u) + \left(\frac{\partial H(\boldsymbol y,u)}{\partial \boldsymbol y} \right)^T \boldsymbol{\bar y}_{k} .
\end{align}
The adjoint iteration \eqref{unsteady_simult_adjoint} can be generated automatically by applying AD to the primal iteration and evaluating the objective function.

\item \textit{Design equation}: The design equation refers to stationarity of the Lagrangian with respect to design changes and reads
\begin{align*}
	0 = \nabla_{u} J_N(\boldsymbol y, u) + \left(\frac{\partial H}{\partial u}\right)^T\boldsymbol{\bar y}.
\end{align*}
It is solved iteratively by
\begin{align}
 \text{iterate } k=0,1\dots \, : \notag \\
	u_{k+1} &= u_k - B_k^{-1} \left(  \nabla_{u}  J_N(\boldsymbol y,u_k) +  \partial_u H(\boldsymbol y,u_k)^T\boldsymbol{\bar y} \right)
\end{align}
with a preconditioning matrix $B_k$.
\end{enumerate}

In the {single-step} one-shot method, the set of all three equations is solved simultaneously in one single coupled iteration:
\begin{align}
	\begin{bmatrix} \boldsymbol y_{k+1} \\ \boldsymbol{\bar y_{k+1}} \\u_{k+1} \end{bmatrix} = \begin{bmatrix} H(\boldsymbol y_k,u_k)  \\ 
				\nabla_{\boldsymbol y} J_N(\boldsymbol y_k,u_k) + \partial_{\boldsymbol y} H(\boldsymbol y_k,u_k) ^T \boldsymbol {\bar y_{k}} \\
				u_k - B_k^{-1} \left(  \nabla_{u}  J_N(\boldsymbol y_k,u_k) +  \partial_u H(\boldsymbol y_k,u_k)^T\boldsymbol{\bar y}_k \right)
	\end{bmatrix}.
\end{align}
As in the steady case, the matrix $B_k$ ensures convergence of the {single-step} one-shot method. In order to approximate the Hessian of the corresponding augmented Lagrangian \eqref{def:augmented_Lagrangian} with respect to design changes, $B_k$ is updated using BFGS-updates on the gradient $\nabla_u L^a$ in each iteration. However, in contrast to {single-step} one-shot optimization with steady state PDEs, the update of the primal and adjoint iteration now each involve a loop over the entire time domain. Only this makes it possible to compute an approximation of the reduced gradient which is used in the design update. 

Enhancing a standard simulation code for {single-step} one-shot optimization involves only minor changes to the time-marching scheme. Since the inner iterator $G^i$ is used in a black box manner,  the stability and robustness properties of the CFD solver are preserved within the new scheme.

\section{{Single-step} one-shot optimization with unsteady RANS}
\label{sec:unsteadyRANS}


The proposed one-shot method is applied to an optimal active flow control problem of unsteady flow around a 2D cylinder. The flow is governed by the unsteady incompressible Reynolds-averaged Navier-Stokes equations (unsteady RANS) {which are solved on a block-structured grid with $12640$ control volumes}. 15 actuation slits are installed on the surface of the cylinder where pulsed actuation is applied according to 
\begin{align}
	 a^l:=u^l\sin {(2\pi ft)} - u^l  \quad \text{for } l=1,\dots,15
\end{align}
while the frequency $f$ is fixed. The amplitudes $u^l$ at each slit are the design variables. The optimization objective is to find optimal actuation parameters $u=(u^1,\dots, u^{15})$ that reduce vorticity downstream the cylinder.  

The optimization problem under consideration is given by
\begin{align}
	\min_{y,u} \,\frac 1T  \int_0^T &C_d(y(t), u) \, \mathrm{d}t \qquad  \text{subject to}  \\
		\begin{split}
			\frac{\partial  v}{\partial t} - \nu \Delta v + (v\cdot \nabla) v + \nabla p &= g  \qquad \qquad \qquad \text{in } \Omega \times (0,T] \\
			\operatorname{div} v &= 0 \qquad \qquad \qquad \text{in } \Omega \times (0,T] \\
			 v(x,t) &= (a^1, \dots, a^{15}) \quad \text{ on } \Gamma \times (0,T] \\
			 v(x,t) &= 0 \qquad \qquad \qquad \text{on } \partial \Omega \setminus \Gamma \times (0,T] \\
			 v(x,0) &= v^0 \qquad \qquad \quad \,\, \text{ in } \Omega
		\end{split}  
\end{align}
where $C_d(t)$ denotes the drag coefficient around the cylinder. The state function $y$ contains the velocity $v(x,t)$ and the pressure $p(x,t)$ in the domain $\Omega$ for a given force field $g$ and viscosity $\nu>0$. The actuation is applied on $15$ slits $\Gamma$ distributed around the cylinder surface.

The governing flow equations are solved with the industry standard CFD simulation code ELAN \cite{holl2012numerical}. ELAN is a second order finite volume code that approximates the transient term with the implicit BDF-2 scheme. The resulting implicit equations are solved one after another using a SIMPLE scheme variant \cite{Peric:2002}. The SIMPLE scheme is a widely used numerical method for solving pressure-linked equations. In that scheme, pressure correction steps to the velocities are performed iteratively until a pseudo-steady state at that time step is reached. Notice, that all numerical effort that makes the simulation code a stable and robust CFD solver is contained inside the pressure-correction iterations while the outer loop shifts the computed solution in time according to the transient approximation. We therefore identify the inner fixed point iterator $G^i$ as one step of the SIMPLE iteration at time step $t_i$. The fixed point iterator $H$, that is used in the proposed one-shot method, performs a loop over all time steps applying one SIMPLE step at each time step.  

The AD tool Tapenade \cite{Tapenade} is applied {to the iteration function $H$ as well as evaluating the discretized objective function. Its reverse mode automatically generates an iterative procedure for solving the adjoint equation and computing the reduced gradient. The design space preconditioner $B_k$ is approximated using BFGS updates based on the reduced gradient, i.e. $\alpha = \beta = 0$. The one-shot iteration stops, when a certain tolerance on the reduced gradient is reached $\|J_u^T + H_u^T\boldsymbol{\bar y}\| \leq \epsilon$, }while in the present test case $\epsilon = 0.001$ was chosen.

\begin{figure}[h] 
	\center
	\includegraphics[width=0.7\textwidth]{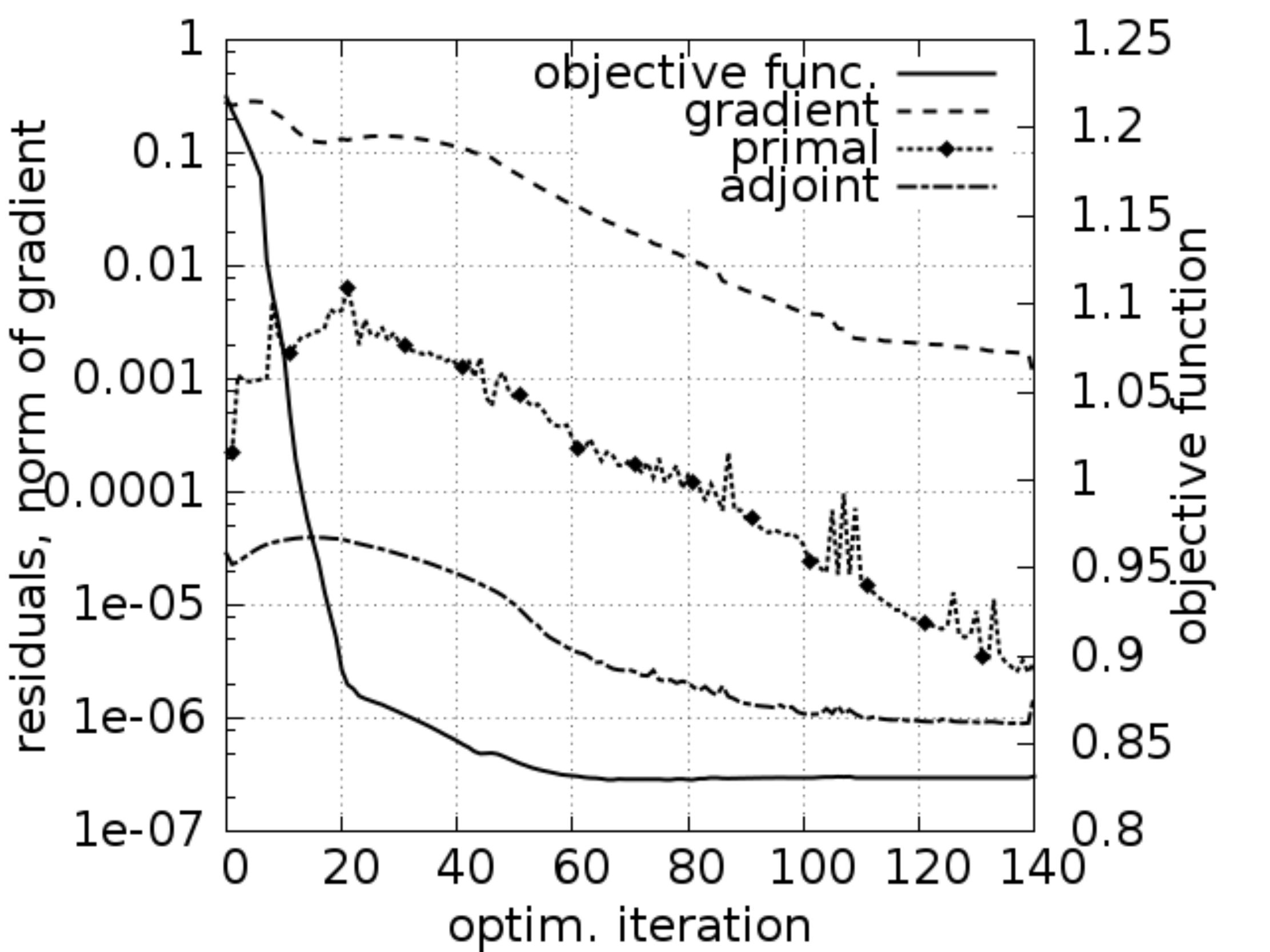}
	\caption{Optimization history of one-shot iterations solving unsteady incompressible Navier-Stokes equations (Re=100)}\label{fig:URANS_optimhistory}
\end{figure}

\begin{figure}[h] 
  \center
  \includegraphics[width=0.7\textwidth]{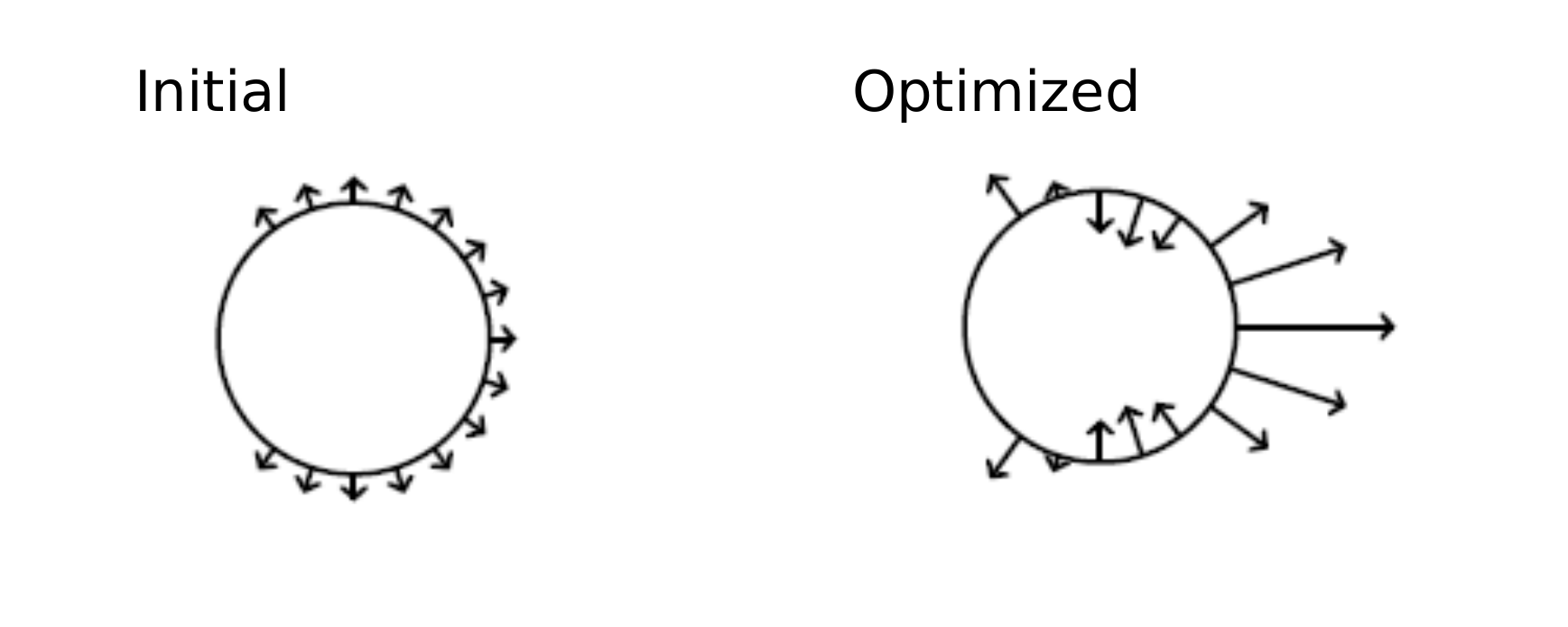}
  \caption{{Initial and optimized amplitudes of the actuation}.}\label{fig:URANS_control}
\end{figure}

Figure \ref{fig:URANS_optimhistory} plots the residuals, the norm of the reduced gradient as well as the objective function during the {single-step} one-shot optimization. An average drag reduction of  about $30\%$ is achieved with the optimization. Primal and adjoint residuals are reduced simultaneously with the reduced gradient indicating the successful application of the {single-step} one-shot method to the unsteady PDE-constrained optimization problem. {The initial and optimized actuation is visualized in Figure \ref{fig:URANS_control}.}

{In this test case, the iteration number needed for convergence of the {single-step} one-shot method was observed to be only 3 times the number of iterations needed for a pure simulation with the modified time-marching scheme with fixed design. This proves the typical bounded retardation of the {single-step} one-shot method as it was observed numerically for the steady case \cite{gauger2009singlestep,gauger2012automated,bosse2014optimal}.}

\section{Improving primal convergence}
\label{sec:improving}

Since the bounded retardation property of the {single-step} one-shot method ensures that the cost for an optimization is only a small factor of the cost of a pure simulation, it is crucial to further investigate the performance of the simulating fixed point solver. We therefore focus on an improvement of the primal iteration by applying adaptive time scales to the state variable.

The fixed point iterations, that solve the unsteady PDE in the {proposed} one-shot method, perform one step of the inner update function $G^i$ at each time step. Due to the inexact approximation of states at previous time steps, the update at a certain time is contaminated by this error. This is also reflected in the lower triangular structure of the Jacobian of $H$. The errors are propagated through the entire time domain and accumulate until the last time step is reached. Thus, the number of iterations needed to reduce the residuals increases as the time domain is enlarged. It was observed numerically, that the dominating contribution to the error occurs in the direction of time while errors in the amplitude of the inexact trajectories are rather marginal. The computed trajectories exhibit a numerical time dilation compared to the final solution approximation.

To reduce the time dilation and improve the primal convergence, we apply an adaptive time scaling approach as introduced in \cite{guenther2014extension}. After each primal update, we assign a trajectory to a scaled time $\tilde t$ such that the new trajectory is in phase with the physical solution. More precisely, we define the new trajectory as
\begin{align}
\tilde y^i := y(\tilde t_i) \quad \forall \, i=1,\dots,N
\end{align} 
where $\tilde t_i$ is chosen such that the residual equation is minimized:
\begin{align}\label{rescaling_minproblem}
	\min_{\tilde t_i} \left\| \frac{y^{i} - y^{i-1}}{\tilde t_{i} - t_{i-1}}  - f_h(y^{i},u) \right\|_2 \quad \forall \, i=1,\dots,N.
\end{align}
The global minimizer of \eqref{rescaling_minproblem} is given by 
\begin{align}\label{rescaling_deftimestep}
	\tilde t_i =  t_{i-1} + \frac{\langle y^i-y^{i-1}, f_h(y^i,u)\rangle}{\|f_h(y^i,u)\|_2^2}.
\end{align}

In this adaptive time scaling approach, we eliminate the error component in the direction of time in such a way, that the numerical time dilation vanishes. {The convergence of the primal iteration is guaranteed with}
\begin{align}
\| \boldsymbol{\tilde y}_{k+1} - \boldsymbol y_* \| &\leq \| \boldsymbol y_{k+1}-\boldsymbol y_*\| \\
 &= \| H(\boldsymbol y_k,u) - \boldsymbol y_*\| \, \overset{k\to \infty}{\longrightarrow} 0
\end{align}
for any design $u\in U$ and the fixed point $\boldsymbol y_*=H(\boldsymbol y_*,u)$.

\section{Numerical results}
\label{sec:numresult}

In this section, we investigate the performance of the fixed point iterator $H$ and the effect of the adaptive time scaling approach. As a first test case, we consider a nonlinear ordinary differential equation (ODE), namely the Van-der-Pol oscillator. Since any unsteady PDE transforms into a system of ODEs after spatial discretization, the Van-der-Pol oscillator is used as a simple model problem. In order to take a step closer to the Navier-Stokes equations, we choose the one-dimensional linear advection-diffusion equation with periodic boundary conditions as a second test case.

\subsection{Van-der-Pol equation}
\label{subsec:vanderpol}
The Van-der-Pol oscillator is a nonlinear oscillator where a damping factor $u \geq 0$ controls the magnitude of the nonlinear term. It can be written as a system of two first order ODEs. With $y=(x,v)^T$ the Van-der-Pol oscillator reads
\begin{align}
	\begin{split}
		 \begin{pmatrix} \dot x(t) \\  \dot v(t)  \end{pmatrix} &= \begin{pmatrix} v(t) \\ -x(t) + u(1-x(t)^2)v(t) \end{pmatrix} \quad \forall t \in (0,T] 	\\
		\begin{pmatrix} x(0)\\v(0) \end{pmatrix} &= \begin{pmatrix} x_0 \\ v_0 \end{pmatrix}
	\end{split}	
\end{align}
where $x$ and $v$ denote the position and the velocity of the oscillator, respectively.

Since we want to resemble the situation where the user is provided with an implicit time-stepping simulation tool, we approximate the transient term with the implicit Backward Euler method. The implicit equations are then solved one after another using an iterative Quasi-Newton method at each time step. According to Section \ref{sec:unsteadyOneshot}, the contractive function $H$ is set up to converge the primal variable simultaneously for all time steps, while $G^i$ represents one step of the Quasi-Newton solver. To remove the numerical time dilation and improve the convergence behavior, the iteration function $H$ is enriched with the adaptive time scaling approach according to \eqref{rescaling_deftimestep}.

\begin{figure}
	\centering
	\includegraphics[width=0.49\textwidth]{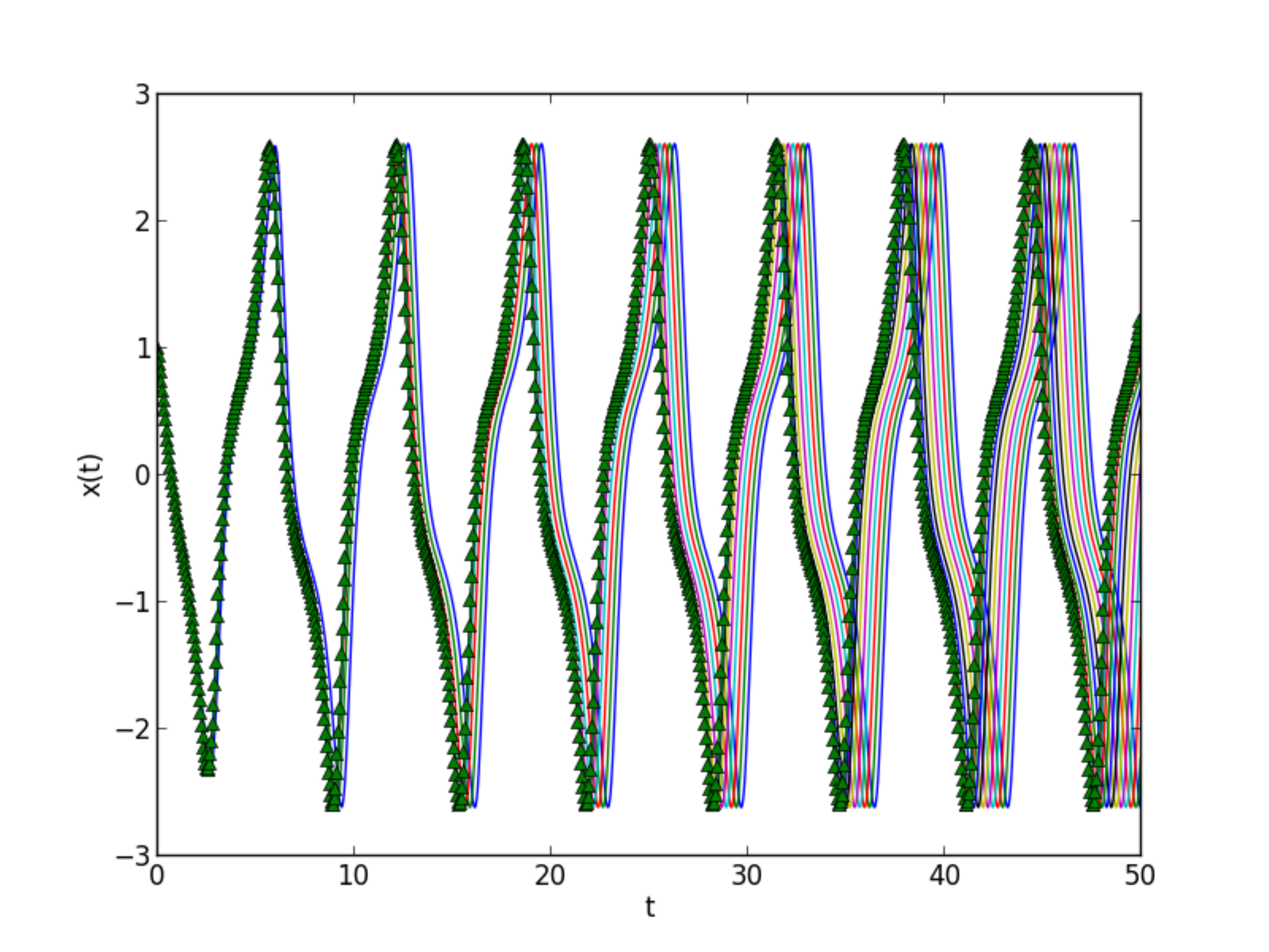}
	\includegraphics[width=0.49\textwidth]{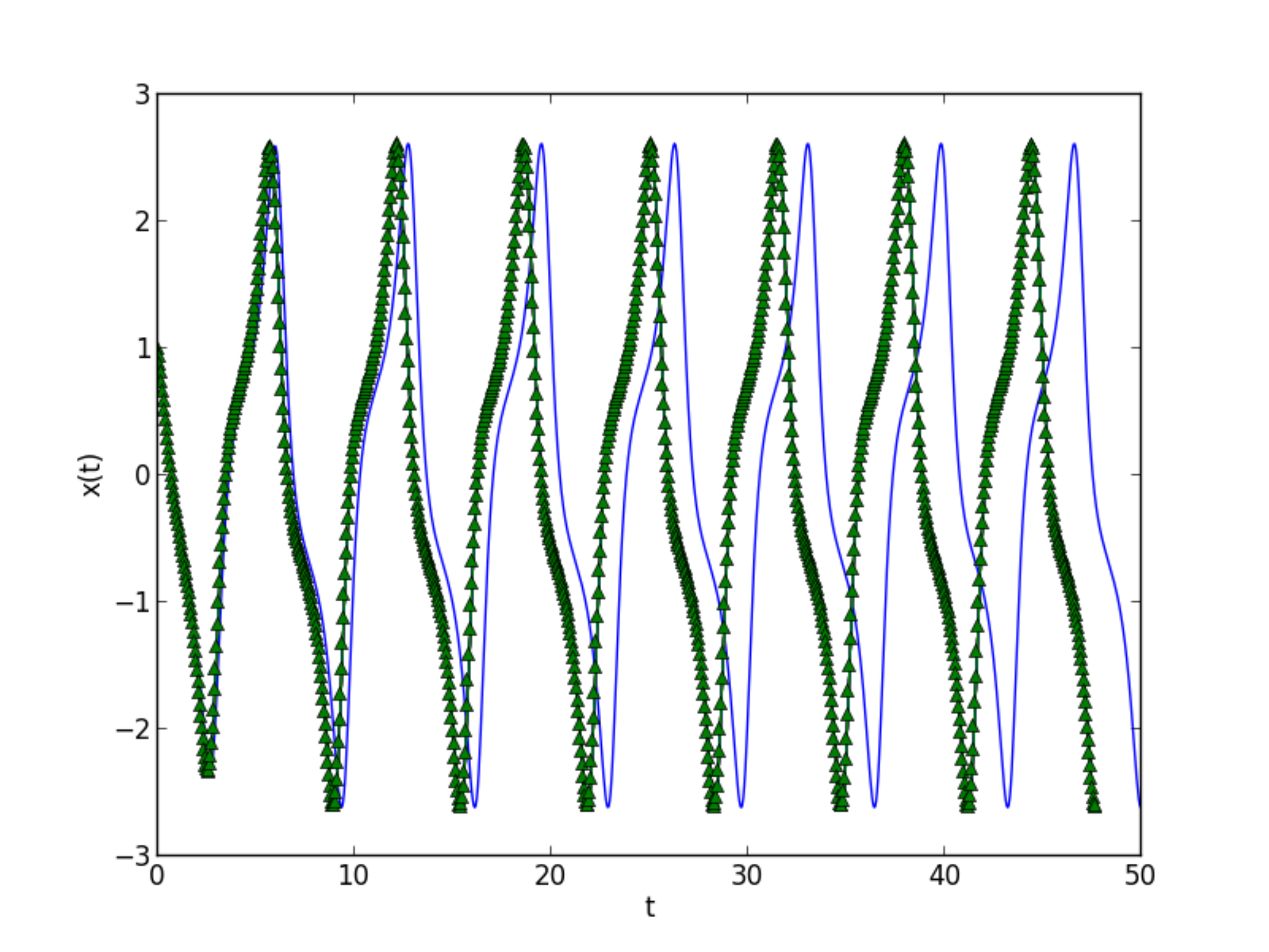}
	\caption{$x$-component of the Van-der-Pol oscillator for $16$ different iterations without (left) and with time scaling approach (right). The physical solution is marked by triangles.}
	\label{fig:VDP_xcomp}
\end{figure}

\begin{figure}
	\centering
	\includegraphics[width=0.49\textwidth]{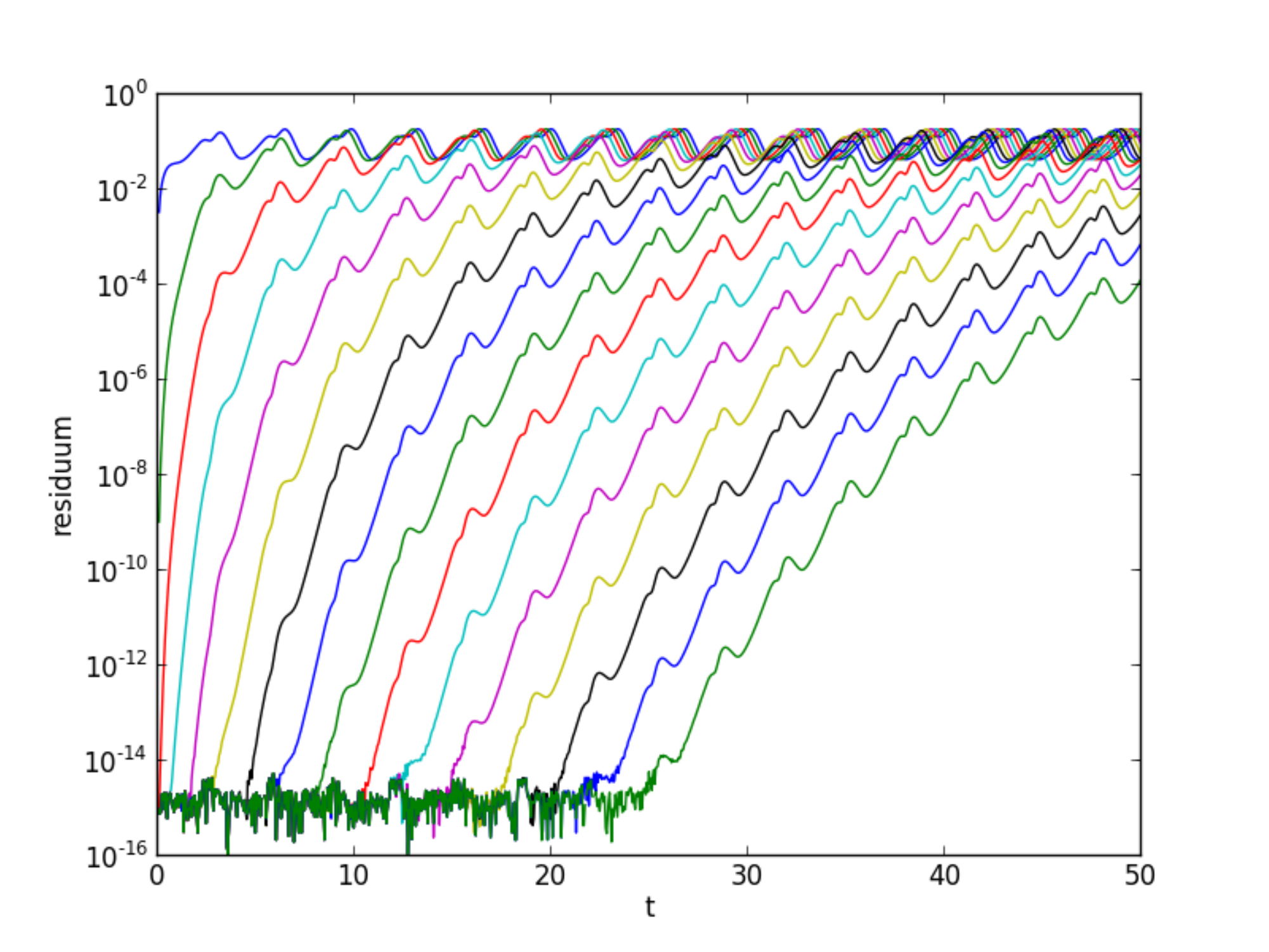} 
	\includegraphics[width=0.49\textwidth]{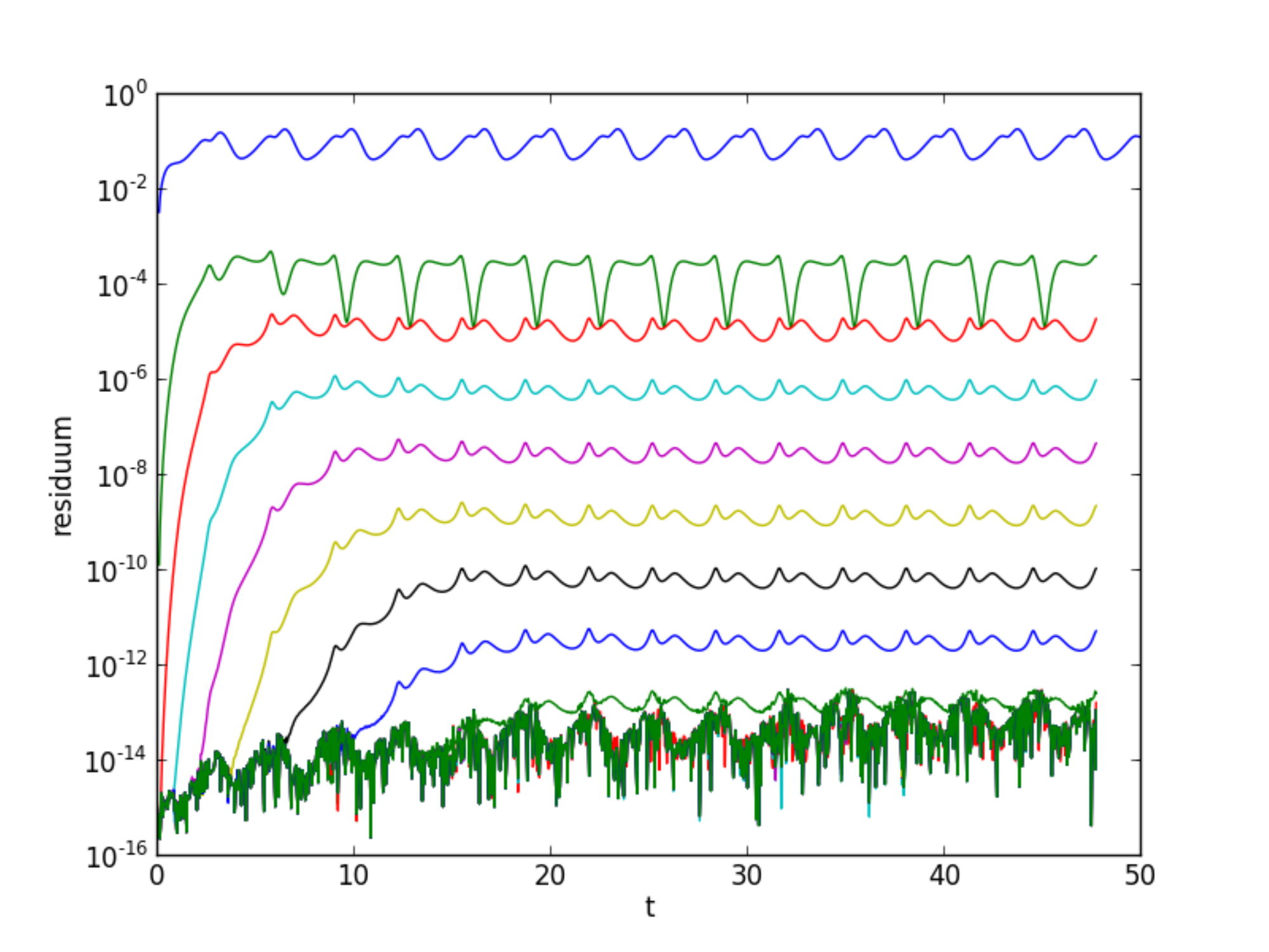}
	\caption{Residuals of the Van-der-Pol oscillator for $16$ different iterations without (left) and with time scaling approach (right).}
	\label{fig:VDP_resid}
\end{figure}

\begin{figure}
	\centering
	\includegraphics[width=0.7\textwidth]{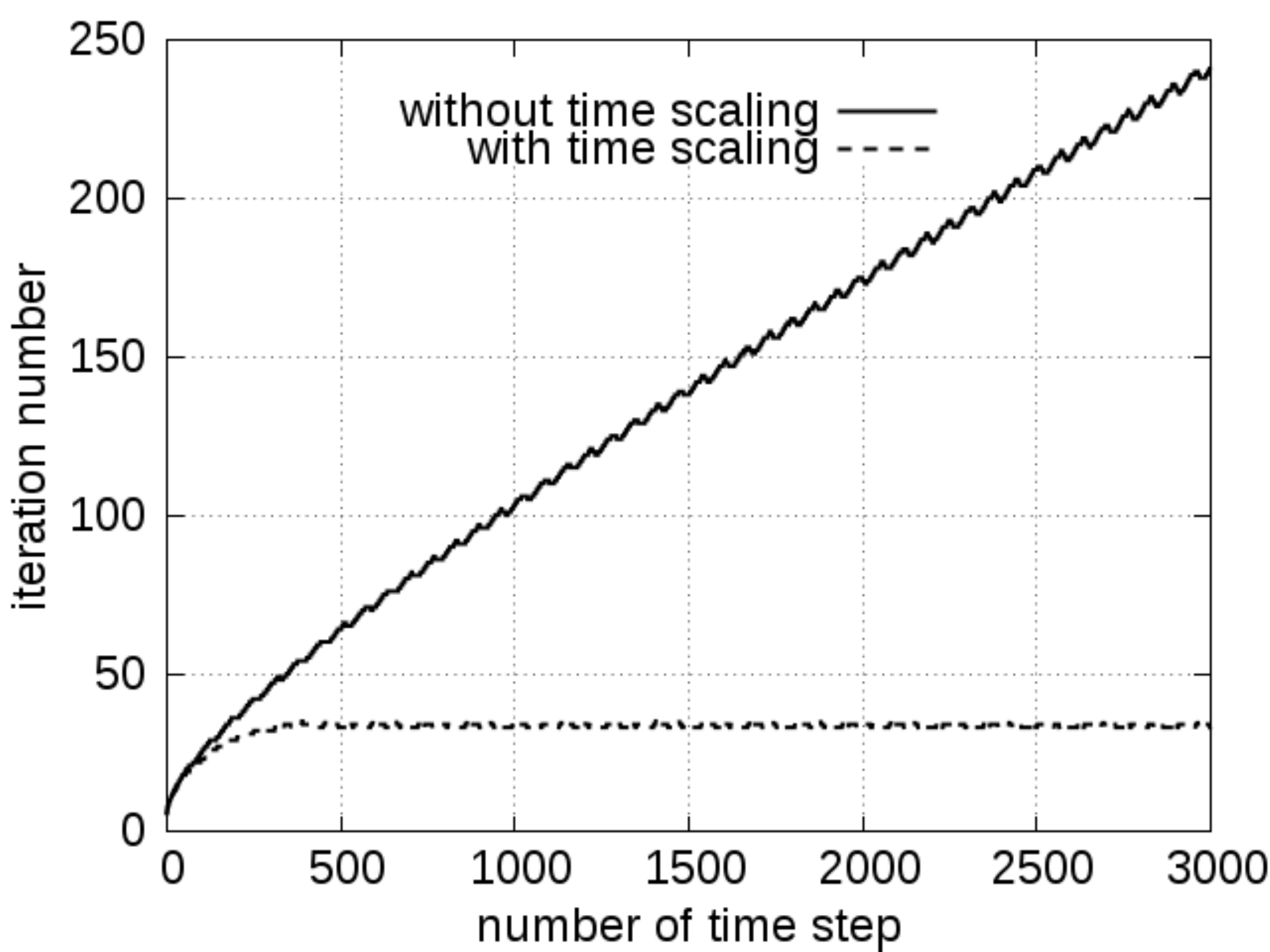}
	\caption{Number of iterations needed for solving Van-der-Pol oscillator with and without time scaling approach.}
	\label{fig:VDP_convergedcompare}
\end{figure}

Figure \ref{fig:VDP_xcomp} visualizes the effect of the time scaling approach on the $x$-component of $16$ intermediate trajectories during the primal computation. On the left, the numerical time dilation can be observed from the spanning bandwidth of intermediate trajectories. In contrast, the time scaled trajectories (right) are all in phase with the physical solution marked by triangles. The corresponding residuals are plotted in Figure \ref{fig:VDP_resid} where the diagonal lines on the left indicate the dependency of the convergence on the time step number. When applying the time scaling approach (right), the residuals drop constantly over the entire time domain. This is also reflected in Figure \ref{fig:VDP_convergedcompare} where the number of iterations that are needed for convergence are plotted for increasing numbers of time steps. Adapting the time scales in each iteration dramatically increases the performance of the primal iteration in this test case.

\subsection{Advection-diffusion equation}
\label{subsec:advectiondiffusion}
As a second test case we investigate the advection-diffusion equation with periodic boundary conditions:
\begin{align}\label{numerics:AD}
\begin{array}{rll}
	 	\partial_t y(t,x)+a \partial_xy(t,x) - \mu \partial_{xx} y(t,x) &=0 & \forall \, x\in (0,1), t\in(0,T]  \\
		 y(0,x) &= h(x)  & \forall \, x\in [0,1]  \\
 		 y(t,0) &= y(t,1) & \forall \, t >0 
 \end{array}
\end{align}
with $a=1$, $\mu = 10^{-5}$ and initial condition $h(x) = \sin(2\pi x)$. We use the Backward Euler discretization in time with initial $\Delta t = 0.01$ and central finite differences in space to approximate the derivative terms. A Quasi-Newton fixed point solver is applied at each time step to solve the implicit equations and resemble the scenario of a general implicit time-marching simulation tool. The iteration function $H$ then performs one update of the Quasi-Newton solver at each time step according to section \ref{sec:unsteadyOneshot}. 

\begin{figure}
	\centering
	\includegraphics[width=0.49\textwidth]{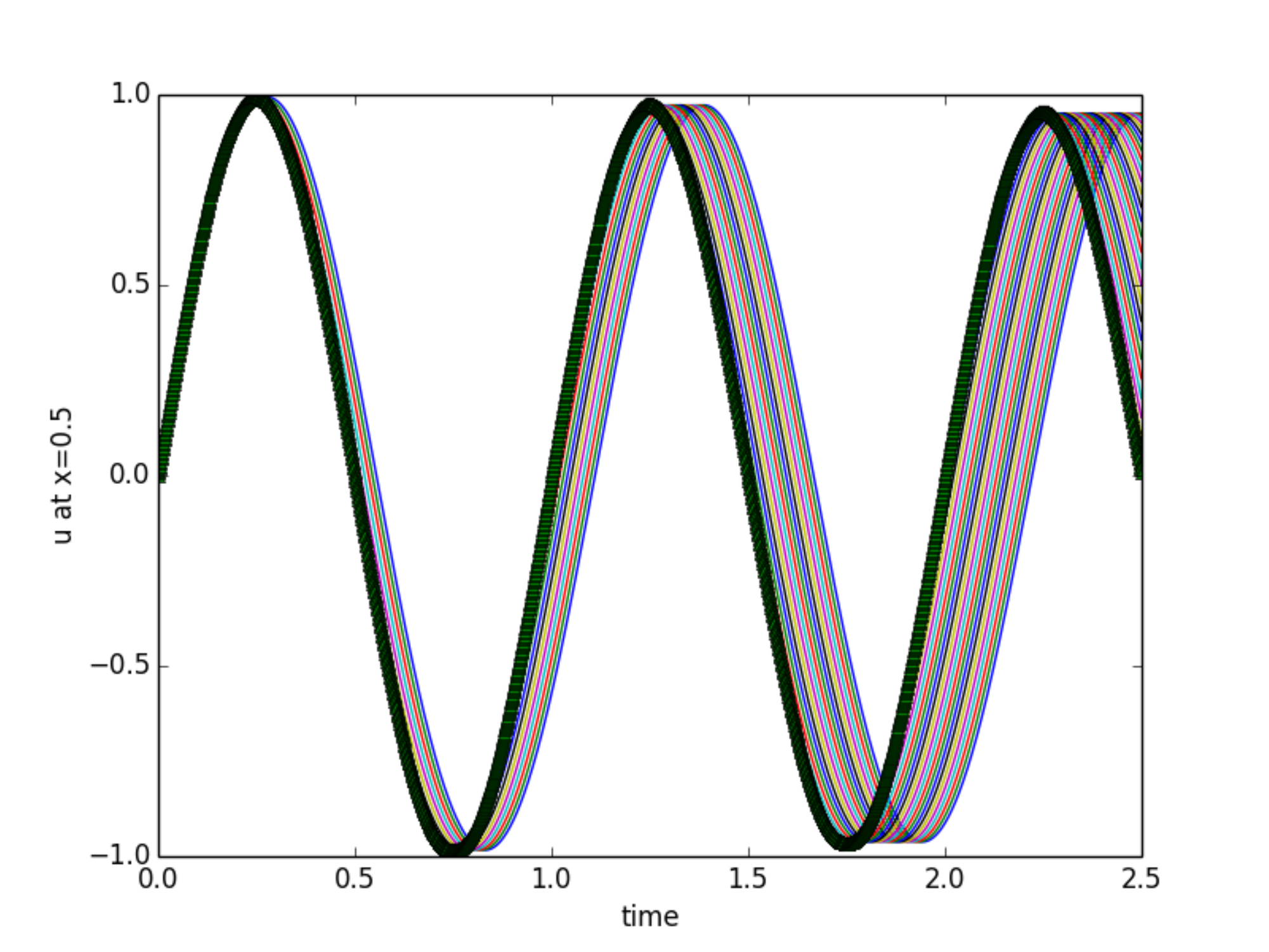} 
	\includegraphics[width=0.49\textwidth]{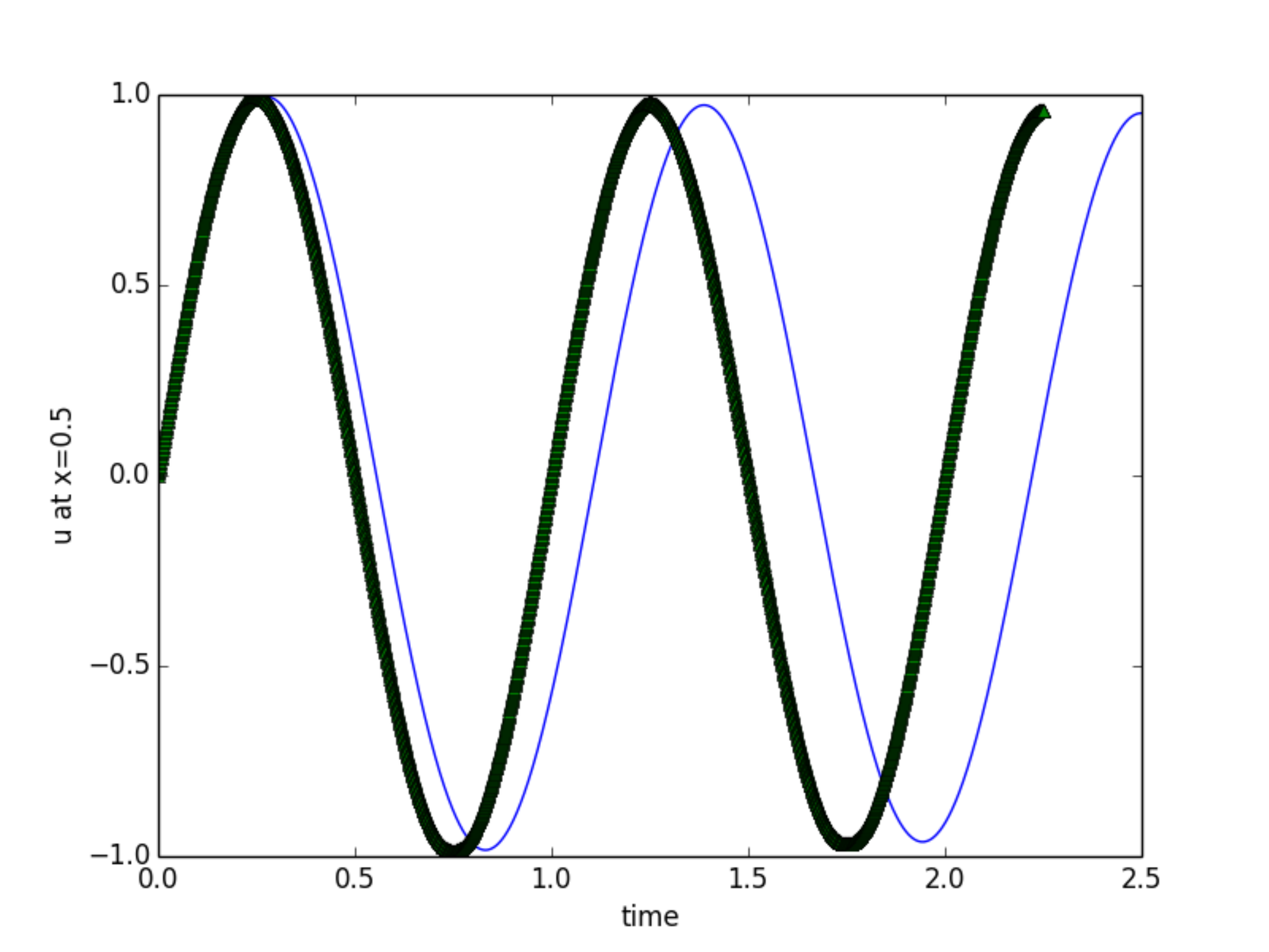}
	\caption{State of the advection-diffusion equation at $x=0.5$ for $16$ different iterations without (left) and with (right) time scaling approach. Physical solution is bold.}	
	\label{fig:AD_state}
\end{figure}

\begin{figure} 
	\centering
	\includegraphics[width=0.49\textwidth]{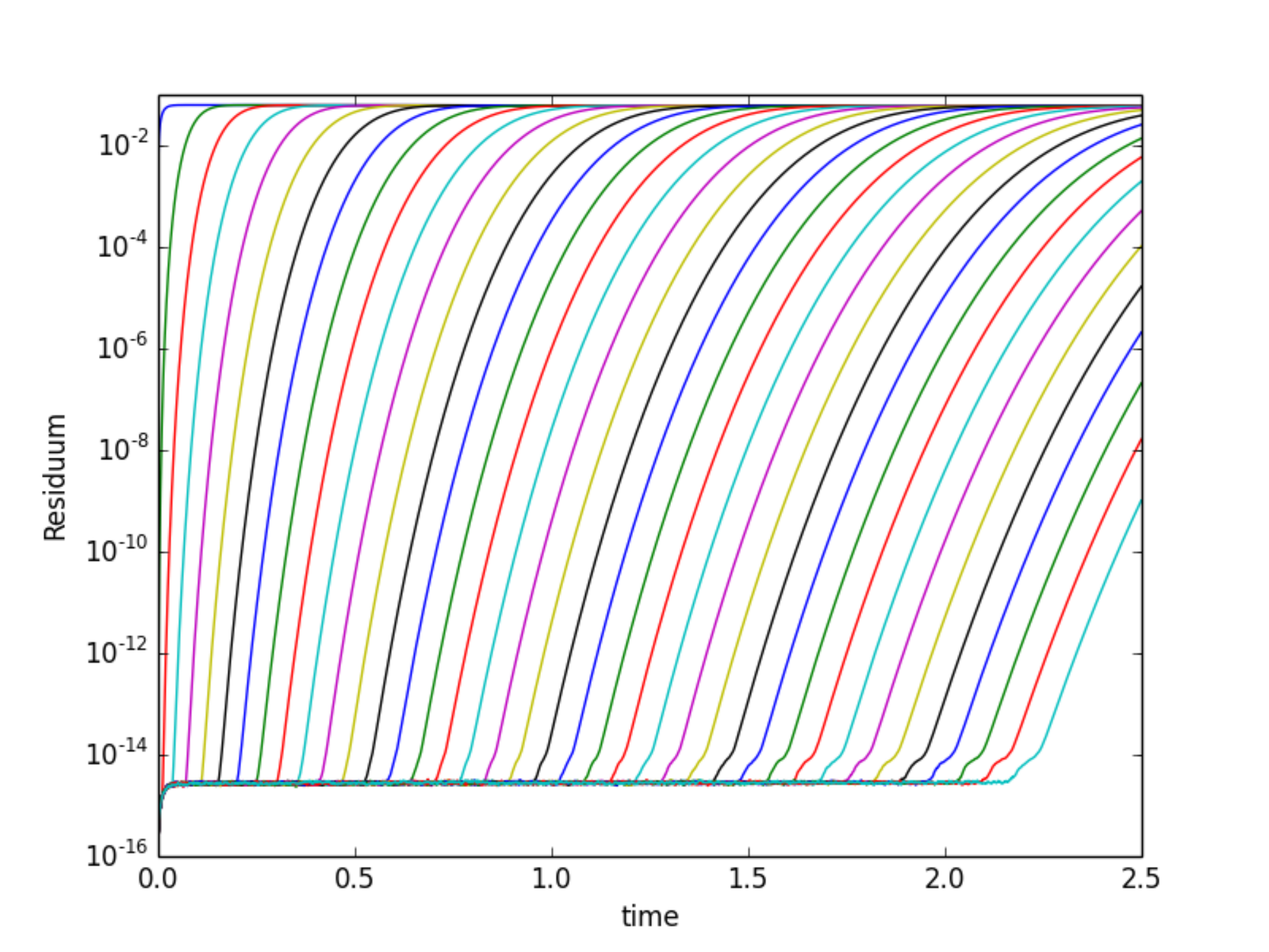} 
	\includegraphics[width=0.49\textwidth]{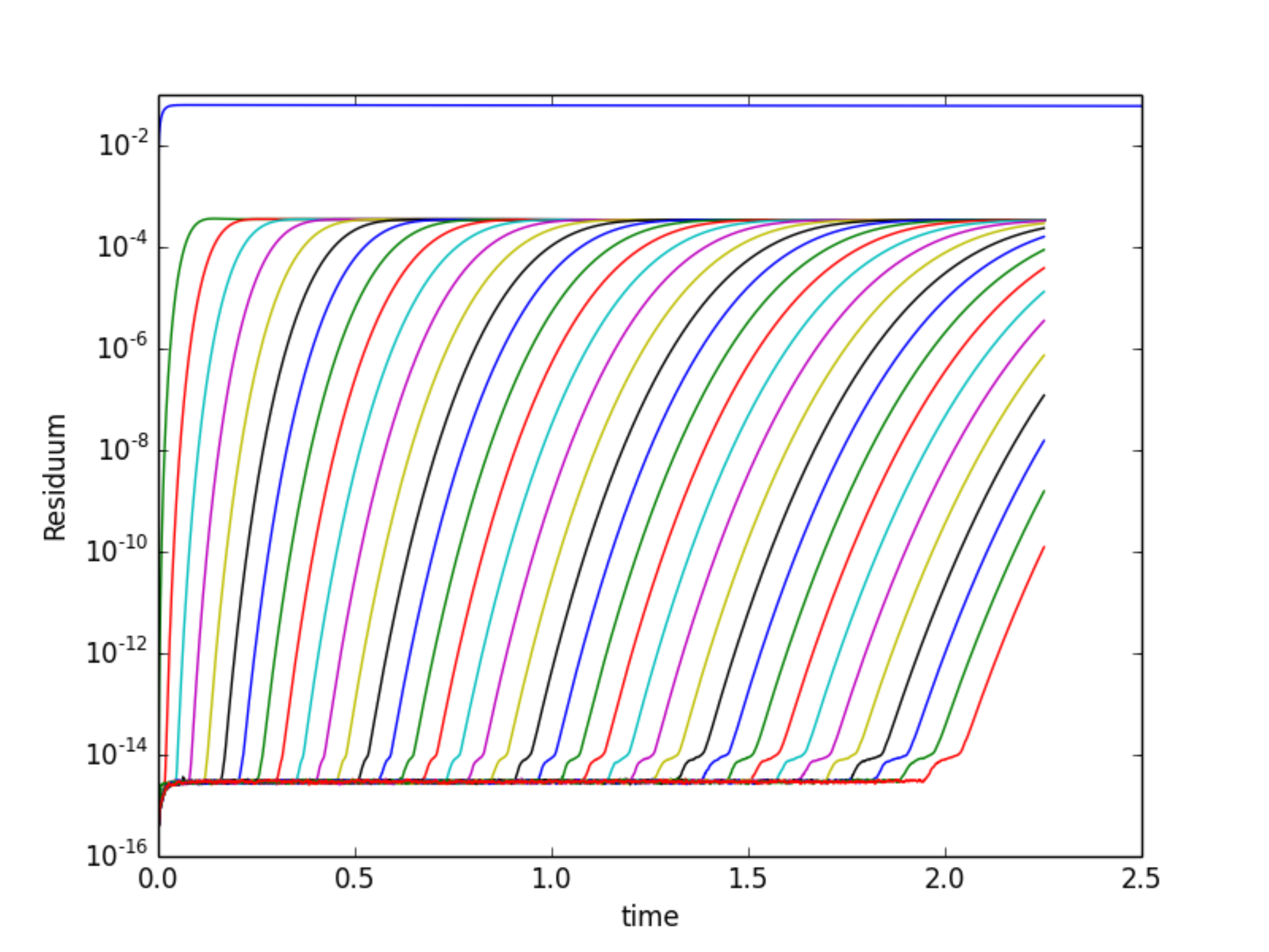} 
	\caption{Residuals of the advection-diffusion equation for $16$ different iterations without (left) and with time scaling approach (right).}
	\label{fig:AD_resid}
\end{figure}

Figure \ref{fig:AD_state} plots $16$ intermediate state trajectories at $x=0.5$ while solving \eqref{numerics:AD} with the fixed point iterator $H$. The numerical time dilation is obvious on the left and enlarges with time. Applying the time scaling approach (right) removes the time dilation such that all intermediate trajectories are in phase with the physical solution. The corresponding residuals are plotted in Figure \ref{fig:AD_resid} where the effect of the time scaling approach is visible on the right: The residuals drop with two orders of magnitude for all time steps after the first time scaling. Nevertheless, the dependency of the convergence on the number of time steps can still be observed from the diagonal structure of the residuals. 

In this test case, applying the time scaling approach still yields an improvement for the primal iteration, but the independence of the iteration number for convergence from the number of time steps, as it was observed in the previous test case, could not be recovered for advection-driven flow. 

\section{Conclusion}\label{sec:conclusion}
In the {single-step} one-shot method, the necessary optimality conditions for PDE-constrained optimization are solved simultaneously in the full space for the primal, the adjoint and the design equations. It is especially tailored for problems where the user is provided with a simulation tool that is to be used in a more or less black box fashion. Simulation tools for steady state PDEs often apply an explicit pseudo-time-stepping scheme which iterates in time until a steady state is reached. In the {single-step} one-shot optimization, these fixed point iterations are enriched by an adjoint and a design update step such that feasibility and optimality is reached simultaneously. The cost of an optimization with the {single-step} one-shot method has proven to be only a small multiple of one pure simulation.

However, if the PDE is unsteady, setting up an appropriate fixed point solver is non-trivial since common schemes often apply an implicit time-marching method and solve the residual equations one after another forward in time. It has been shown in this paper, that these time-marching schemes can be modified to fit into the proposed one-shot optimization framework by reducing the number of inner iterations for solving the implicit equations. In the resulting approach, the entire time trajectory of the unsteady PDE is updated within one iteration. Applying AD to the modified time-marching scheme as well as evaluating the discrete objective function automatically generates a consistent discrete adjoint iteration. Augmenting the modified primal iteration with an adjoint update and a preconditioned reduced gradient step for the design yields the {single-step} one-shot optimization method for unsteady PDE-constrained optimization. The method has been applied to an optimal active flow control problem governed by the unsteady RANS equations. 

The modified time-marching scheme has been further improved applying adaptive time scales. In that approach, the intermediate trajectories are shifted in time such that they are in phase with the physical solution. For a general model problem, the time scaling approach yields an improved convergence behavior that is independent of the time domain resolution. An application to unsteady flow driven by advection showed that the time scaling approach still yields an improvement for the convergence speed while the dependence on the number of time steps remains.

%



  \bibliographystyle{elsarticle-num} 
 \bibliography{mybib}





\end{document}